\documentclass[a4paper]{article}
\pdfoutput=1

\usepackage{IEK10} 
\usepackage{natbib}

\makeatother
\usepackage{amssymb}
\usepackage{xcolor}
\usepackage{mathtools, cuted}
\usepackage{tabularx}
\usepackage{eurosym}
\usepackage[utf8]{inputenc}

\usepackage{pstricks, pst-plot}	
\usepackage{graphicx} 
\usepackage{wrapfig} 
\usepackage[figuresright]{rotating}
\usepackage{lipsum}
\usepackage{placeins}
\usepackage{etoolbox}

\usepackage{amsmath}
\numberwithin{equation}{section} 
\usepackage{amssymb}
\usepackage{setspace}
\usepackage{enumerate}
\usepackage{enumitem}
\usepackage[english]{babel}
\usepackage{blindtext}
\usepackage{epsfig} 
\usepackage[utf8]{inputenc}			
\usepackage{ifpdf} 							
\usepackage{multicol}						
\usepackage{titletoc}						
\usepackage{setspace}						
\usepackage{ulem} 							
\usepackage{longtable}						
\usepackage{multirow,bigdelim}
\usepackage{booktabs} 						
\usepackage{parcolumns}						
\usepackage{tabularx} 						
\usepackage{threeparttable}			  
\usepackage{graphicx}							
\usepackage{flafter}							
\usepackage{placeins}							
\usepackage{rotating}							
\usepackage{tikz}								
\usetikzlibrary{external}										
\usepackage{pgfplots}							
\usepackage{epstopdf}							
\pgfplotsset{compat=newest}						
\pgfplotsset{plot coordinates/math parser=false}
\usepackage{pdfpages}							
\graphicspath{{../02_Grafiken/}} 
\usepackage{float}
\usepackage{mhchem}

\renewcommand{\min}[1][]{
	\ifthenelse{\isempty{#1}}{\operatorname{min}}{\ensuremath{\underset{#1}{\text{min}\,}}}
}

\def\LTT{
  Institute of Technical Thermodynamics, 
  RWTH Aachen University,
  Schinkelstraße 8,
  52062 Aachen,
  Germany
}

\newcommand{\mytitle}{Gotta catch ‘em all: Modeling All Discrete Alternatives for Industrial Energy System Transitions}

\newcommand{\affil}{
  \begin{itemize}[leftmargin=3mm, itemsep=0mm]
    \item[$^a$]\LTT
  \end{itemize}
}

\def\firstAuthor{Hendrik Schricker}
\newcommand{\myauthor}{\firstAuthor$^{a}$, Benedikt Schuler$^{a}$, Christiane Reinert$^{a}$, and Niklas von der Aßen$^{a,*}$}
\newcommand{\correspondingauthor}{N. von der Aßen, \LTT \newline E-mail: niklas.vonderassen@ltt.rwth-aachen.de}
\author{\myauthor}

\usepackage[
  colorlinks,
  linkcolor=blue,
  citecolor=blue,
  urlcolor=blue,
  pdftitle={\mytitle},
  pdfauthor={\firstAuthor}
]{hyperref}
\usepackage[capitalise, nameinlink]{cleveref}
\crefname{table}{Tab.}{Tab.}

\newcommand{\setpgfexternalcounter}[1]{
  \makeatletter%
  \pgfkeysgetvalue{/tikz/external/figure name}\myexternalname
  \expandafter\gdef\csname c@tikzext@no@\myexternalname\endcsname{#1}%
  \makeatother
}

\begin{document}

  \thispagestyle{firststyle}

  \begin{center}
    \begin{large}
      \textbf{\mytitle}
    \end{large} \\
    \myauthor
  \end{center}

  \vspace{0.5cm}

  \begin{footnotesize}
    \affil
  \end{footnotesize}

  \vspace{0.5cm}

  \begin{abstract}
Industrial decision-makers often base decisions on mathematical optimization models to achieve cost-efficient design solutions in energy transitions. However, since a model can only approximate reality, the optimal solution is not necessarily the best real-world energy system.
Exploring near-optimal design spaces, e.g., by the Modeling All Alternatives (MAA) method, provides a more holistic view of decision alternatives beyond the cost-optimal solution. However, the MAA method misses out on discrete in-vestment decisions. Incorporating such discrete investment decisions is crucial when modeling industrial energy systems.   

Our work extends the MAA method by integrating discrete design decisions. We optimize the design and operation of an industrial energy system transformation using a mixed-integer linear program. First, we explore the continuous, near-optimal design space by applying the MAA method. Thereafter, we sample all discrete design alternatives from the continuous, near-optimal design space. 

In a case study, we apply our method to identify all near-optimal design alternatives of an industrial energy system. We find 128 near-optimal design alternatives where costs are allowed to increase to a maximum of one percent offering decision-makers more flexibility in their investment decisions. Our work enables the analysis of discrete design alternatives for industrial energy transitions and supports the decision-making process for investments in energy infrastructure.
\end{abstract}

\vspace{0.5cm}


\vspace{0.75cm}

\newpage
\pagestyle{laterstyle}
\section{Planning industrial energy system transformations}\label{sec:Intro}
Industrial processes are typically energy-intensive. Often, the energy demand is provided by on-site, decentralized energy systems. Decarbonizing industrial energy supply is critical in global greenhouse gas emission reduction strategies. However, planning the optimal transition of industrial energy systems is challenging. Energy system optimization models (ESOMs) have been applied to guide decision-makers in recent years. 

However, the calculated optimum might not always be the most suitable solution, as optimization models represent only an approximation of the real world and rely on many uncertain assumptions for input parameters (parametric uncertainty) and for model equations (structural uncertainty). Hence, the appropriate handling of uncertainties has become a major challenge for future ESOMs \citep{Fodstad.2022}. 

For handling parametric uncertainties, several approaches are suggested in the literature \citep{Yue.2018}. However, few methods exist for incorporating structural uncertainty. A method to deal with structural uncertainty in ESOMs is Modeling to Generate Alternatives (MGA) \citep{DeCarolis.2017,Voll.2015}. The MGA method increases insights in decision flexibility for decision-makers by exploring alternative solutions in the near-optimal solution space. Nevertheless, the MGA method does not cover the near-optimal solution space evenly.  Furthermore, the generated alternatives are limited in their number and often represent extreme alternatives near the edge of the near-optimal solution space. \cite{Pedersen.2021} propose the Modeling All Alternatives (MAA) method to overcome these disadvantages. The MAA method explicitly calculates a geometric representation of the near-optimal solution space and, thus, enables uniform coverage. However, accurate modeling of industrial energy systems requires the inclusion of integer variables since most supply technologies are only available in discrete capacity steps. 

In this work, we identify all discrete, near-optimal design alternatives. Our method provides decision-makers with a set of suitable alternatives. To the best of our knowledge, our work is the first extension of the MAA method for a sector-coupled model, which can identify all discrete design alternatives for industrial energy systems. 

\FloatBarrier
\section[Modeling all discrete design alternatives]{Modeling all discrete design alternatives}\label{sec:Method}
In our work, we sequentially solve linear energy system design optimization problems to determine near-optimal design alternatives regarding the minimal total annualized cost $TAC$. The energy system has to cover temporally resolved model-exogenous energy demands. As variables, we consider the capacity expansion $d\in\mathcal{D}$ of potentially built energy conversion or storage components and their respective operation $o \in \mathcal{O}$.

In this work, we restrict the capacity expansion $d_c \in \mathcal{D}_c$ for each component $c \in \mathcal{C}$ to discrete capacity steps $d_c^\mathrm{min}$. The discrete set of installable capacities is thus defined as 
\begin{align}
\mathcal{D}_c = \{d_c=\lambda \cdot d_c^\mathrm{min}: \lambda \in \mathbb{N}_0 \} \; \forall c \in \mathcal{C}.
\end{align}
Overall, the cost-optimal design is determined by solving the mixed-integer linear program 
\begin{align} TAC^*=\min_{d\in\mathcal{D}, o \in \mathcal{O}}TAC(d,o).\end{align}
We identify the near-optimal design space $\mathcal{W}_\varepsilon$, which we define as 
\begin{align}
\mathcal{W}_\varepsilon=\{d\in\mathcal{D} \mid \min_{o\in\mathcal{O}}TAC(d,o)\le TAC^*\cdot(1+\varepsilon)\}.
\end{align} 

The MAA method builds upon the assumption that the considered optimization problem and thus, the near-optimal solution space $\mathcal{W}_\varepsilon$ is convex. However, the domain set of the design variables $\mathcal{D}$ is non-convex per definition. Therefore, we relax the discrete character of $\mathcal{D}$ by introducing the relaxed domain set of the design variables $\widetilde{\mathcal{D}} = \mathbb{R}_{\ge 0}^{\vert\mathcal{C}\vert}$. We obtain the continuous, near-optimal design space $\widetilde{\mathcal{W}}_\varepsilon$ to which we can apply the MAA method (cf. Section \ref{sec:contSpace}). We then reintroduce the discrete character of the installed capacities $d$ and employ an efficient sampling method to identify all discrete design alternatives (cf. Section \ref{sec:discrSpace}).  

\subsection{Exploring the continuous, near-optimal design space}\label{sec:contSpace}
To span the continuous, near-optimal design space $\widetilde{\mathcal{D}}$, we iteratively solve the following optimization problem with altered search directions $n\in\mathbb{R}^{\vert\mathcal{C}\vert}$:
\begin{align} 
v(n)=\max_{d\in\widetilde{\mathcal{W}}_\varepsilon}n^\mathrm{T}d. 
\end{align}
We obtain a vertex $v\in\mathcal{V}$ of the continuous, near-optimal design space $\widetilde{\mathcal{W}}_\varepsilon$ for each search direction $n$. In the first iterations, we choose the positive and negative unit vectors as search direction $n$, i.e., we minimize and maximize the installed capacity $d_c$ of each component $c\in\mathcal{C}$. To speed up the computation time, we parallelize the exploration of the search directions on 48 cores. The normal directions of the facets of the convex hull $\mathcal{L}$ of all so far found vertices serve as search directions for the later iterations. 
We introduce a new termination criterion to quantify the degree of exploration of the continuous, near-optimal solution space: We define a hyperplane for each search direction $n$ and vertex $v(n)$, which is normal to $n$ and contains $v(n)$. The combination of all these hyperplanes defines a polyhedron $\mathcal{U}$. The volume $V$ of the near-optimal space $\widetilde{\mathcal{W}}_\varepsilon$ is limited by the volume of $\mathcal{U}$ as an upper bound and the volume of $\mathcal{L}$ as a lower bound. We stop exploring search directions if  
\begin{align} 
\frac{V(\mathcal{U})-V(\mathcal{L})}{V(\mathcal{U})}\le \delta
\end{align} with a user-defined tolerance $\delta$ holds true, i.e., if at least a ratio of $1-\delta$ of the volume of $\widetilde{\mathcal{W}}_\varepsilon$ has been explored. With our termination criterion, we enable quantifying the degree of exploration of near-optimal solution spaces for the first time.

\subsection[Modeling all discrete design alternatives]{Discretizing the near-optimal design space}\label{sec:discrSpace}
From the convex hull $\mathcal{L}$ of vertices of the continuous, near-optimal design space $\widetilde{\mathcal{W}}_\varepsilon$, we derive the defining set of linear equations for the discrete, near-optimal design space  
\begin{align} 
\mathcal{W}_\varepsilon=\{d \in \mathcal{D} | Ad\le b\}.
\end{align}  
Instead of validating the feasibility of all discrete designs $d\in\mathcal{D}$ in a brute force approach, we implement a recursive sorting algorithm that efficiently identifies all discrete, near-optimal designs: First, we determine the smallest hyperrectangle $\mathcal{H}$ with vertices $h\in\mathcal{D}$ that contains $\widetilde{\mathcal{W}}_\varepsilon$. We then split $\mathcal{H}$ into smaller hyperrectangles $\mathcal{H}_i, i=1 \dots m$ and validate the feasibility of the discrete design alternatives corresponding to the vertices of these smaller hyperrectangles. If $Ah \le b$ holds true for all vertices $h$ of a hyperrectangle $\mathcal{H}_i$, all discrete designs $d\in\mathcal{H}_i$ are near-optimal due to convexity of $\widetilde{\mathcal{W}}_\varepsilon$ and $\mathcal{H}_i$. If none of the vertices $h$ satisfies $Ah \le b$, we exclude all discrete designs $d \in \mathcal{H}_i$ from further consideration. If at least one vertex $h$ but not all vertices satisfy $Ah \le b$, we again split $\mathcal{H}_i$ into smaller hyperrectangles until all discrete designs $d \in \mathcal{H}$ are sorted. 

The recursive approach has a slightly lower accuracy than the brute force approach, since some rare edge cases are excluded, reducing the identified near-optimal solution space by a very small extent. However, we drastically decrease the required computation time and RAM especially for high-dimensional and highly resolved solution spaces. 

\newpage
\FloatBarrier
\section[Case Study]{Case Study: Transformation of an Industrial Energy System under Structural Uncertainty}\label{sec:casestudy}
\begin{figure}[ht]
	\centering
	\includegraphics[width=1.0\textwidth]{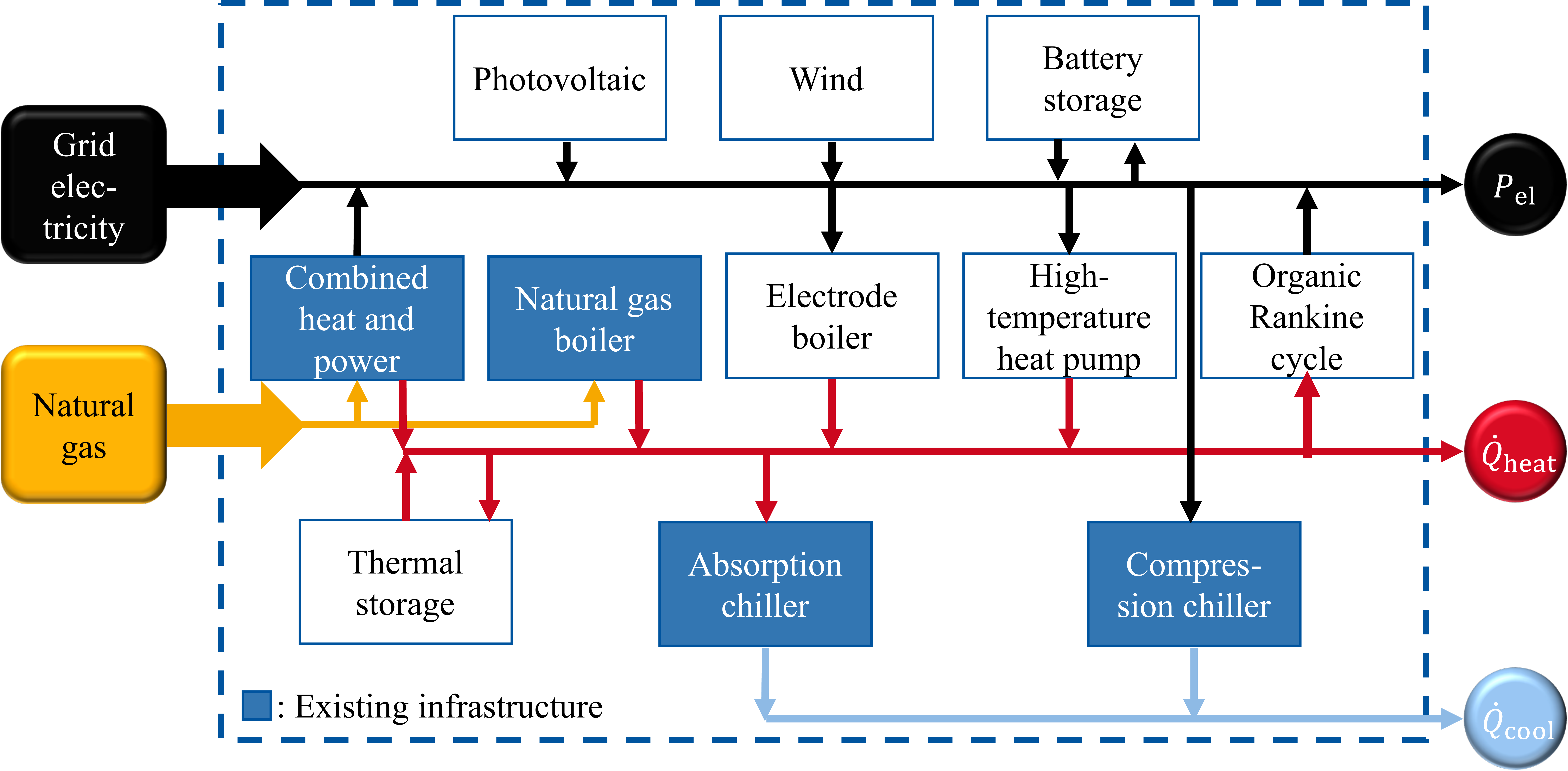}
	\caption{Structure of the industrial energy system providing electricity, heating, and cooling energy \citep{Voll.2013,Baumgartner.2019,Reinert.2023}. Blue boxes are existing infrastructure, and white boxes are technology investment options.}	
	\label{fig:case_study}
\end{figure}

The industrial energy system analyzed in this study (Figure \ref{fig:case_study}) was first published by \cite{Voll.2013} and extended by \cite{Baumgartner.2019} and \cite{Reinert.2023}. We identify all near-optimal design alternatives for the transformation of an existing industrial energy system. We consider the total annualized cost $TAC$ including \ce{CO2} certificate cost as objective function. The energy system serves temporally resolved electricity, heating, and cooling demands. Overproduction of electricity, heating, and cooling is not allowed. Natural gas and electricity can be imported from the grid. We assume temporally resolved electricity and natural gas prices and a constant price for \ce{CO2} certificates. 
The existing energy system infrastructure is mainly fossil-fuel based. Installable technologies are photovoltaic, wind turbines, sector-coupling and storage technologies. All technologies have underlying techno-economic data, e.g., specific investment costs and technical efficiencies, and environmental impact data. We assume fixed capacity steps within the range of typical sizes for each installable technology. Furthermore, we do not limit the installable capacities by a maximum installable capacity.
We implement our method with the energy system optimization framework SecMOD \citep{Reinert.2022}. In each optimization, we consider one year of operation, aggregated into four typical days with eight typical hours each. We assume a relative objective slack of $\varepsilon=1~\%$ and evaluate 2144 search directions until we cover at least 95~\% of the volume of the continuous, near-optimal design space. 

\begin{figure}[ht]
	\centering
	\includegraphics[width=1.0\textwidth]{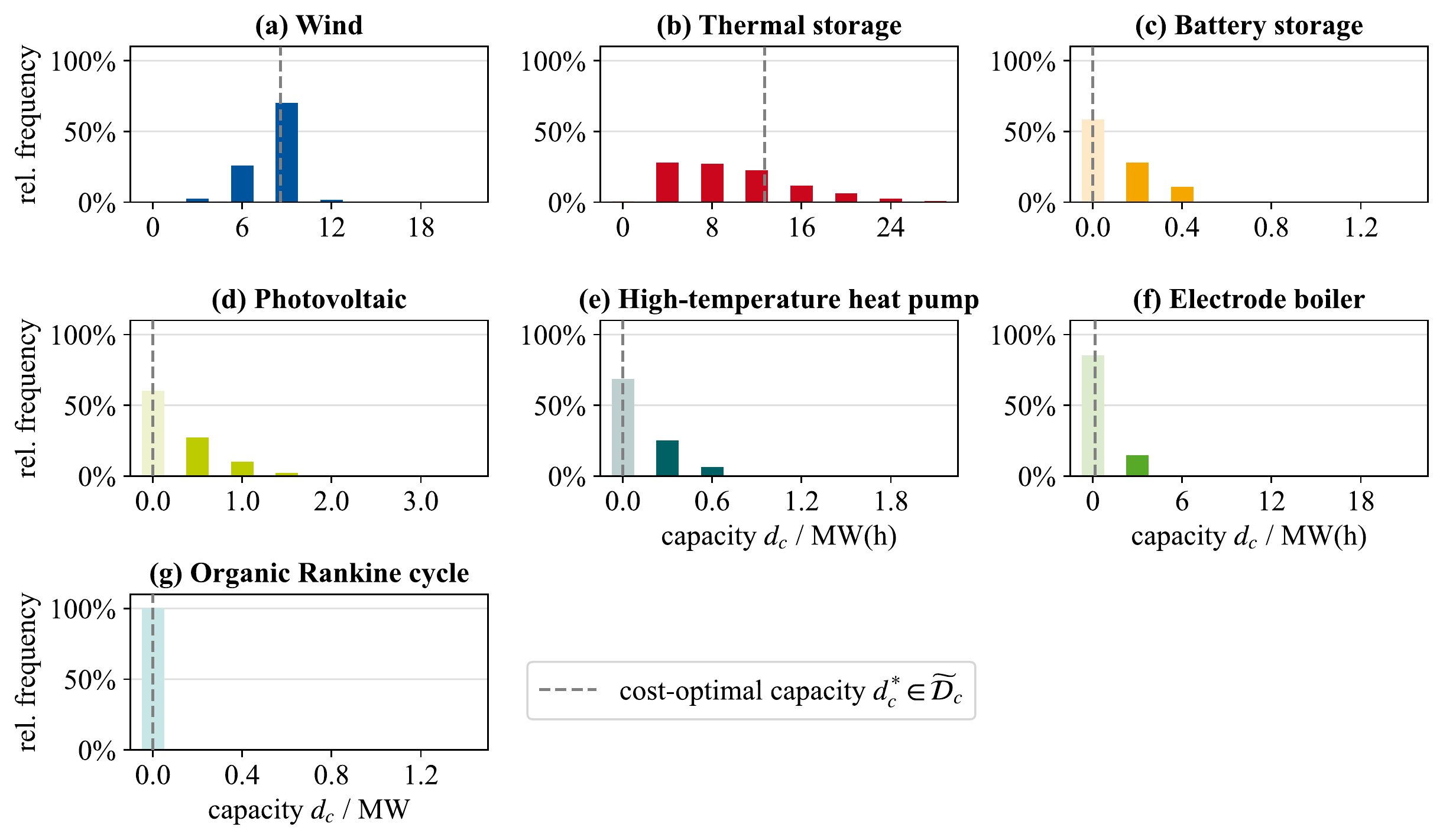}
	\caption{Relative frequency of installed capacity of each technology across all identified near-optimal design alternatives. The dotted vertical line indicates the cost-optimal solution for the relaxed problem with continuous capacities.}	
	\label{fig:rel_freq}
\end{figure}

In total, we identify 128 near-optimal design alternatives. We derive the frequency distribution of the installed capacity per technology across all near-optimal design alternatives to analyze the resulting capacity ranges for each technology in comparison to the cost-optimal solution from the relaxed problem (cf. Figure \ref{fig:rel_freq}). In the cost-optimal design, wind turbines and a thermal storage are installed. The thermal storage enables temporally limited overproduction of heat from the CHP plants and thus allows for electricity-demand-driven operation of the CHP plants. Both investments thus, lead to a reduction of costly electricity imports. 
The frequency distributions of the wind and thermal storage capacity reveal that the installation of wind turbines and the thermal storage are necessary conditions for near-optimality in the given case study. The ORC is not installed in any of the 128 design alternatives. Thus, all potential designs that include an ORC are at least 1~\% more costly. The remaining technologies can be installed as backup capacity or for converting excess electricity into heat, which can be used for demand coverage or be stored in the thermal storage. However, these investments are not necessary to achieve near-optimality regarding the minimal total annualized cost but offer decision-makers more flexibility.

\FloatBarrier	
\section{Conclusions}\label{sec:conclusion}
In this work, we present an extension of the Modeling All Alternatives method for discrete technology investment options. In the first step, we systematically explore the near-optimal design to obtain an approximated geometric representation of the continuous, near-optimal design space according to the original Modeling All Alternatives method. We introduce an improved termination criterion for the exploration of the continuous, near-optimal design space, which enables to quantify the degree of exploration of near-optimal solution spaces for the first time. To account for discrete investment options in the second step, we implement a computationally efficient sampling method and thus enable complete enumeration of all discrete design alter-natives. 
Our work allows incorporating structural model uncertainty for industrial energy system design optimizations. Furthermore, we identify decision flexibilities for near-optimal investments thus support decision-makers in transforming their energy system.

\section*{Acknowledgements}
HS gratefully acknowledges the financial support of the Ministry of Economics, Industry, Climate Protection and Energy of North-Rhine Westphalia (Grant number: EFO 0133E). Simulations were performed with computing resources granted by RWTH Aachen University under project ID 6307.

\FloatBarrier

  \bibliographystyle{apalike}
  \renewcommand{\refname}{Bibliography}  
  \bibliography{literature.bib}

\begin{thebibliography}{}

\bibitem[Baumg{\"a}rtner et~al., 2019]{Baumgartner.2019}
Baumg{\"a}rtner, N., Delorme, R., Hennen, M., and Bardow, A. (2019).
\newblock Design of low-carbon utility systems: Exploiting time-dependent grid
  emissions for climate-friendly demand-side management.
\newblock {\em Applied Energy}, 247:755--765.

\bibitem[DeCarolis et~al., 2017]{DeCarolis.2017}
DeCarolis, J., Daly, H., Dodds, P., Keppo, I., Li, F., McDowall, W., Pye, S.,
  Strachan, N., Trutnevyte, E., Usher, W., Winning, M., Yeh, S., and Zeyringer,
  M. (2017).
\newblock Formalizing best practice for energy system optimization modelling.
\newblock {\em Applied Energy}, 194:184--198.

\bibitem[Fodstad et~al., 2022]{Fodstad.2022}
Fodstad, M., {Del Crespo Granado}, P., Hellemo, L., Knudsen, B.~R., Pisciella,
  P., Silvast, A., Bordin, C., Schmidt, S., and Straus, J. (2022).
\newblock Next frontiers in energy system modelling: A review on challenges and
  the state of the art.
\newblock {\em Renewable and Sustainable Energy Reviews}, 160:112246.

\bibitem[Pedersen et~al., 2021]{Pedersen.2021}
Pedersen, T.~T., Victoria, M., Rasmussen, M.~G., and Andresen, G.~B. (2021).
\newblock Modeling all alternative solutions for highly renewable energy
  systems.
\newblock {\em Energy}, 234:121294.

\bibitem[Reinert et~al., 2023]{Reinert.2023}
Reinert, C., Nolzen, N., Frohmann, J., Tillmanns, D., and Bardow, A. (2023).
\newblock Design of low-carbon multi-energy systems in the secmod framework by
  combining milp optimization and life-cycle assessment.
\newblock {\em Computers {\&} Chemical Engineering}, 172:108176.

\bibitem[Reinert et~al., 2022]{Reinert.2022}
Reinert, C., Schellhas, L., Mannhardt, J., Shu, D.~Y., K{\"a}mper, A.,
  Baumg{\"a}rtner, N., Deutz, S., and Bardow, A. (2022).
\newblock Secmod: An open-source modular framework combining multi-sector
  system optimization and life-cycle assessment.
\newblock {\em Frontiers in Energy Research}, 10:637.

\bibitem[Voll et~al., 2015]{Voll.2015}
Voll, P., Jennings, M., Hennen, M., Shah, N., and Bardow, A. (2015).
\newblock The optimum is not enough: A near-optimal solution paradigm for
  energy systems synthesis.
\newblock {\em Energy}, 82:446--456.

\bibitem[Voll et~al., 2013]{Voll.2013}
Voll, P., Klaffke, C., Hennen, M., and Bardow, A. (2013).
\newblock Automated superstructure-based synthesis and optimization of
  distributed energy supply systems.
\newblock {\em Energy}, 50(C):374--388.

\bibitem[Yue et~al., 2018]{Yue.2018}
Yue, X., Pye, S., DeCarolis, J., Li, F.~G., Rogan, F., and Gallach{\'o}ir,
  B.~{\'O}. (2018).
\newblock A review of approaches to uncertainty assessment in energy system
  optimization models.
\newblock {\em Energy Strategy Reviews}, 21:204--217.

\end{thebibliography}

\end{document}